\documentstyle[12pt]{amsart}
\makeatletter
\def\LaTeX{\leavevmode L\raise.42ex
    \hbox{\kern-.3em\size{\sf@size}{0pt}\selectfont A}\kern-.15em\TeX}
\makeatother

\newcommand{\BibTeX}{{\rm B\kern-.05em{\sc i\kern-.025emb}\kern-.08em\TeX}}

\newtheorem{proposition}{Proposition}[section]
\newtheorem{thm}[proposition]{Theorem}
\newtheorem{defn}[proposition]{Definition}
\newtheorem{corollary}[proposition]{Corollary}
\newtheorem{question}[proposition]{Question}

\newtheorem{result}[proposition]{Result}
\newtheorem{problem}[proposition]{Problem}
\newsymbol\square 1003
\newcommand{\N}{\Bbb{N}}
\newcommand{\rnk}{\mbox{rnk}}
\renewcommand{\span}{\mbox{span}}
\font\sll=cmti12

\makeatletter
\def\@currentlabel{2.1}\label{e:dispaa}
\def\@currentlabel{2.21}\label{e:dispau}
\def\@currentlabel{2.22}\label{e:dispav}
\def\@currentlabel{2.23}\label{e:dispaw}
\def\@currentlabel{2.24}\label{e:dispax}
\def\theequation{\thesection.\@arabic\c@equation}

\makeatother

\makeatletter
\def\alphenumi{%
  \def\theenumi{\alph{enumi}}%
  \def\p@enumi{\theenumi}%
  \def\labelenumi{(\@alph\c@enumi)}}
\makeatother

\begin{document}

\title[The Homogeneous Banach Space Problem]%
{Some Questions Arising from the Homogeneous Banach Space Problem}
\author{Peter G. Casazza}
\address{Department of Mathematics, University of Missouri,
Columbia, Missouri 65211}
\email{mathpc@@umcvmb.bitnet}
%  \thanks will become a 1st page footnote.
%  Don't type a period at the end; it will be supplied.
\thanks{Supported by NSF DMS 9001344 and a grant from the University of
Missouri Research Council}
\thanks{Part of this research was conducted while the author was a
visiting scholar at the University of Cambridge, ENGLAND}
%\author{Author Two}
%\address{Mathematical Research Section, School of Mathematical Sciences,
% Australian National University, Canberra ACT 2601, Australia}
%\email{ABC@@mathsci.anu.edu.au}
%\thanks{The final version of this paper will be submitted for
% publication elsewhere}

%\keywords{Author's key words go here}
%  Math Subject Classifications

\maketitle

\begin{abstract}
We review the current state of the homogeneous Banach space problem.  We
then formulate several questions which arise naturally from this problem,
some of which seem to be fundamental but new.  We give many examples
defining the bounds on the problem.  We end with a simple construction
showing that every infinite dimensional Banach space contains a subspace
on which weak properties have become stable (under passing to further
subspaces).  Implications of this construction are considered.
\end{abstract}

\section{Introduction}

A Banach space is said to be \underline{homogeneous} if it is isomorphic
to all of its infinite dimensional subspaces.

\bigskip
\noindent\underline{\bf The Homogeneous Banach Space (HBS) Problem}
\bigskip

Is every homogeneous Banach space isomorphic to a Hilbert space?

\bigskip
This problem has frequently been referred to as ``Banach's problem''
since it is stated in Banach's book [B].  However, it was pointed
out in [MTJ2] that in the original polish version of his book,
Banach attributes this problem to Mazur.  The objective of this paper is
to generate renewed interest in the HBS problem and a host of other
interesting problems surrounding it.

Essentially no progress was made on this problem of mazur until recently
when J.~Bourgain [B$^*$] solved the finite dimensional version and
W. B. Johnson [J1] solved a special case of the general problem.
Both papers rely on advances made in the local theory of Banach spaces in
the 1980's.  These results are discussed in section 3.

In section2, we consider the long list of immediate questions we can also
ask about a homogeneous Banach space, most of which are still unanswered.
In section 4, we will look at a series of problems which are stronger
than the homogeneous Banach space problem, in the sense that a positive
answer to one of these would yield a positive answer to the HBS problem.
In section 5 we look at some fundamental questions which arise from the
HBS problem, some of which do not seem to have been asked before.
Finally, in section 6, we give a simple construction that can be carried
out in any Banach space to yield a subspace which is stable for weak
properties, for passing to further subspaces.  We will also show that
this is the strongest subspace that we can expect to find inside of every
Banach space.

The author expresses his deepest gratitude to W. B. Johnson for many
interesting discussions related to the material in this paper.

\section{Some Immediate Questions}

There are several immediate questions which we can ask about a
homogeneous Banach space $X$: (1) Is $X$ reflexive?; (2) Is $X$
suprreflexive?; (3) Does $X$ have a separable dual space?  None of these
questions has been answered yet.  However, we do not believe that any of
these questions will be important for the solution of the problem.  But
these questions do give rise to some interesting problems in the general
theory of Banach spaces which are discussed in section 5.  Even the
isometric version of the HBS problem is still open.

\bigskip
\noindent\underline{\bf The Homogeneous Banach Space (HBS) Problem
(Isometric Version)}
\bigskip

If $X$ is isometric to all of its infinite dimensional subspaces, is $X$
isometric (or even isomorphic) to a Hilbert space?

\bigskip
With so little progress having been made on this problem, we might hope
to get some movement on it by strengthening our hypotheses.
Unfortunately, assuming that $X$ is homogeneous and has a unconditional
basis, or even a symmetric basis, does not seem to help.  These
hypotheses do yield that $X$ is superreflexive and that $X$ is weak
cotype 2 [J1] but do not seem to lead to any serious breakthrough
on the problem.

Another obvious question is: Is every homogeneous Banach space $X$
isomorphic to a square?  (i.e. Is $X \approx X \oplus X$?)  Even such an
elementary question was unresolved until recently and requires some heavy
machinery of W. T. Gowers and B. Maurey [GoM].  Recall that a
projection
$P$ on a Banach space is said to be \underline{nontrivial} if $\rnk
P=\rnk
(I-P)=\infty$.

\begin{thm} [G${}_0$M]
If an infinite dimensional Banach space $X$ has no nontrivial projections
on any infinite dimensional subspace, then every bounded operator on $X$
is a strictly singular pertubation of the identity.
\label{thm:2.1}
\end{thm}

This means that every bounded operator $T$ on such a space can be written
as: $T=aI+S$, where $S$ is a strictly singular operator (i.e. $S$ is not
an isomorphism when restricted to any infinite dimensional subspace of
$X$.)  We also need the notion of Fredholm index.  If $T:X\to Y$ is a
bounded operator with closed range, put $\alpha(T)=\dim\ker T$,
$\beta(T)=\dim Y/TX$.
If either $\alpha(T) < \infty$ or $\beta(T) < \infty$, we define the
\underline{Fredholm index} $i(T)$ by: $i(T)=\alpha(T)-\beta(T)$.  If
$i(T)$ is defined and is finite, then $T$ is called a \underline{Fredholm
operator}.  (See chapter 2.c of \cite{LT1}, for basic information on
$i(T)$).  We are now ready for:

\begin{result}\label{re:2.2}
If $X$ is homogeneous then $X \approx X \oplus X$.
\end{result}

\begin{pf}
If not, then $X$ has no non--trivial projection on any infinite
dimensional subspace.  Let $T$ be an isomorphism of $X$ onto a hyperplane
in $X$.  By Theorem~\ref{thm:2.1}, $T=aI+S$.  But, $T$ has Fredholm index
$-1$ while $aI+S$ has Fredholm index $0$ by Proposition 2.c.10 of [LT1].
This contradiction completes the proof.
\end{pf}

There are some immediate questions which would be useful for the solution
to the HBS problem.  One such question is,

\begin{question}\label{q:2.3}
If $X$ is homogeneous, is $X^*$ homogeneous?
\end{question}

In section 3, we will see a case where W. B. Johnson \cite{J1} made use
of such a hypothesis.  If $X$ is homogeneous then, since $X$ has a
subspace with a basis, every subspace of $X$ has a basis.  In particular,
$X$ is separable and every subspace of $X$ has the approximation
property.  It follows (\cite{LT2}, Theorem 1.g.6) that;

\begin{equation}
\sup\{p|X \mbox{ is type } p\}=2=\inf\{q|X \mbox{ is cotype }q\}.
\label{eq:2.4}
\end{equation}

This fact gives rise to some interesting open questions in section 4.

We also feel that the following open question will be important to the
solution of the HBS problem.

\begin{question}
If $X$ is homogeneous, is $X$ uniformly isomorphic to all of its infinite
dimensional subspaces?
\label{q:2.5}
\end{question}

That is, does there exist a $K\geq 1$ so that $X$ is $K$--isomorphic to
all of its infinite dimensional subspaces?  One reason for the importance
of question~\ref{q:2.5} is that, perhaps the HBS problem has a positive
answer
in the uniform case but a negative answer in the general case.  We can
prove, with some relatively soft infinite dimensional theory, the
following result:

\begin{proposition}
If $X$ is a homogeneous Banach space, then there is a constant $K\geq 1$
so that $X$ $K$--embeds into every infinite dimensional subspace of $X$.
\label{p:2.6}
\end{proposition}

Proposition~\ref{p:2.6} will follow from a general result on minimal
Banach
spaces.  A Banach space is \underline{minimal} if it embeds into every
one of its infinite dimensional subspaces.  We discuss such spaces in
more detail in section 6.

\begin{proposition}
If $X$ is a minimal Banach space then there is a $K\geq 1$ so that $X$
$K$--embeds into every infinite dimensional subspace of $X$.
\label{p:2.7}
\end{proposition}

\begin{pf}
Assume, for the sake of contradiction, that $X$ is not uniformly
embeddable into all of its infinite dimensional subspaces.  Then no
subspace of $X$ has this property either.  So there are infinite
dimensional subspaces $Y_1 \supset Y_2 \supset \cdots$ so that $X$ does
not $2^{f(n)}$--embed into $Y_n$, where $f(n)=n(1+\sqrt{n})^2$.  Choose
$y_n\in Y_n$ so that $(y_n)$ is a basic sequence in $X$ and let
$E_n=\span_{n\leq i}\, y_i$.  By assumption, there is a subspace
$Z\subset E_1$ and a $K\geq 1$ so that $Z$ is $K$--isomorphic to $X$.
Without loss of generality, we may assume that $\dim\, E_1 / Z=\infty$.
For each $n$ let $H_n=Z\cap E_n$ and choose $F_n \subset E_n$ with $H_n
\subset F_n$ and $k=\dim E_n/F_n=\dim Z/H_n\leq n$.  There are
projections $P:E_n \to F_n$ and $Q:Z\to H_n$ with $\|p\|\leq 1+\sqrt{n}$
and $\|Q\|\leq 1+\sqrt{n}$.  Hence, there are $k$--dimensional spaces
$W\subset E_n$ and $W'\subset Z$ so that $W'\oplus_1 H_n$ is
$(1+\sqrt{n})$--isomorphic to $Z$ while $span (W,\, H_n)$ is
$(1+\sqrt{n})$--isomorphic to $W\oplus_1\, H_n$.  Since $W$ and $W'$ are
$n$--isomorphic, it follows that $V=\span (W,\, H_n)$ is
$n(1+\sqrt{n})^2$--isomorphic to $Z$ and hence $V$ is $K n
(1+\sqrt{n})^2$--isomorphic to $X$.  But $V\subset E_n$ implies $X$ is $K
n (1+\sqrt{n})^2$--embeddable into $E_n$.  For large $n$, this
contradicts our assumption that $X$ is not $2^{f(n)}$--embeddable into
$E_n$, for $f(n)=n(1+\sqrt{n})^2$.
\end{pf}

Proposition~\ref{p:2.7} easily gives a subspace $Y\subset X$ with a
finite dimensional decomposition $\Sigma\oplus E_n$ so that for every
$k$,
$(E_n)^{\infty}_{n=k}$ is dense in the family of all finite dimensional
subspaces of $X$.  Hence, there is a constant $K\geq 1$ so that every
finite dimensional subspace of $X$ is $K$--isomorphic to a
$K$--complemented subspace of every finite codimensional subspace of $X$.

Another consequence of a positive answer to 2.4 comes from Krivine's
Theorem \cite{Kr}, \cite{MS}.  If $X$ is homogeneous, a positive answer
to~\ref{q:2.5}, combined with~\ref{eq:2.4} and Krivine's Theorem implies:

\begin{thm}
If $X$ is $K$--isomorphic to all of its infinite dimensional subspaces,
then for every $n$, there is a basis $(x_i)$ of $X$ with basis constant
$\leq K$ and for every finite set of natural numbers $I$ with $|I|=n$,
$(x_i)_{i\in I} \approx_{K}(e^2_i)_{i\in I}$, where $(e_i^2)$ is the unit
vector basis of $l_2$.
\end{thm}

This says that $X$ has a sequence of bases $(x_n^k)_{n=1}^\infty$ with
uniformly bounded basis constants so that every $k$ elements of
$(x_n^k)^{\infty}_{n=1}$ are uniformly the $l_{2}^{k}$ unit vector basis.
This property itself is very strong and we might hope that it already
characterizes a Hilbert space.  Although, as we will see, it does not
characterize Hilbert space, it does give rise to the notion of sequences
of successively better bases for a Banach space.

\begin{defn}
We say that a Banach space $X$ has \underline{basis property $p$} $(1
\leq
p\leq \infty)$ if there is a $K\geq 1$ so that for every $n$, there is a
basis $(x_i)$ of $X$ with basis constant $\leq K$ and for every
$n$--element subset $I$ of the natural numbers, $(x_i)_{i\in I}
\approx_{K}(e_i^p)_{i\in I}$, where $(e^p_i)$ is the unit vector basis of
$l_p$ (or $c_0$ if $p=\infty$).  To simplify notation, we write
$(x_n)\leq (y_n)$, for two sequences in Banach spaces $X$, $Y$
respectively, if there is a constant $K\geq 1$ so that for every sequence
of scalars $(a_n)$, $$\left\|\sum_{n=1}^\infty\, a_n\, x_n \right\| \leq
K \left\|\sum_{n=1}^\infty \, a_n \, y_n\right\|.$$  We now have,
\label{de:2.9}
\end{defn}

\begin{proposition}
Let $X$ be a Banach space with a normalized basis $(x_n)$ satisfying
$(x_n)<(e_n^p)$.  Then $X\oplus l_p$ has basis property $p$.
\label{p:2.10}
\end{proposition}

\begin{pf}
Let $(x_n)$ be the normalized basis of $X$ satisfying $(x_n)<(e_n^p)$.
Fix a natural number $m$.  Define $y_n\in X\oplus l_p$ by:
\begin{eqnarray*}
y_{k2m+(k+1)} &= &\left(x_{k+1}, \dfrac1{(2m)^{1/p}}\,
\sum_{i=k2m+1}^{(k+1)2m}\, e^p_i \right)\\
y_{k2m+(k+1)}+j &= &\left(0, e^p_{k2m+j} \right), \hspace{1in} \mbox{for
}
1\leq j
\leq 2m,
\end{eqnarray*}
and $k=0,1,2,\ldots$.  For any sequence of scalars $(a_i)$ and any $t<s$,
choose $k_1, k_2, j_1, j_2$ so that $t=k_1\, 2m+(k_1+1)+j_1$ and
$s=k_2\, 2m+(k_2+1)+j_2$.  Then, treating $X\oplus l_p$ as an
$l_p$--sum and letting $b$ be the basis constant of $(x_n)$ we
have,
\begin{eqnarray*}
\lefteqn{\left\|\sum_{i=1}^t \, a_i\, y_i \right\|=} \\
& & \left[ \left\| \sum_{i=0}^{k_1}\, a_{i2m+(i+1)}\, x_{i+1} \right\|^p
\right. +  \\
& & \sum_{i=0}^{k_1}\, \sum_{j=1}^{2m}\, \left| a_{i2m+(i+1)+j} -
\dfrac{1}{(2m)^{1/p}}\, a_{i2m+(i+1)} \right|^p  \\
& & +\sum_{j=1}^{j_1} \left| a_{i2m+(i+1)+j}-\dfrac{1}{(2m)^{1/p}}\,
a_{i2m+(i+1)} \right|^p + \\
& & \left. \dfrac{2m-j_1}{2m}\, \left|a_{k_1 2m+(k_1+1)} \right|^p
\right]^{1/p} \\
& & \leq \left[ b \left\| \sum_{i=0}^{k_1}\, a_{i2m+(i+1)}\, x_{i+1}
\right\|^p \right. +  \\
& & \sum_{i=0}^{k_1}\, \sum_{j=1}^{2m} \, \left|
a_{i2m+(i+1)+j}-\dfrac{1}{(2m)^{1/p}}\, a_{i2m+(i+1)}\right|^p \\
& & +\left. \sum_{j=1}^{j_1} \, \left|
a_{i2m+(i+1)+j}-\dfrac{1}{(2m)^{1/p}}\,a_{i2m+(i+1)} \right| \right]^{1/p} \\
& & \leq \left[ b^2 \left\| \sum_{i=0}^{k_1}\, a_{i2m+(i+1)}\, x_{i+1}
\right\|^p \right. + \\
& & \sum_{i=0}^{k_2}\, \sum_{j=1}^{2m}\, \left|
a_{i2m+(i+1)+j}-\dfrac{1}{(2m)^{1/p}}\, a_{i2m+(i+1)} \right|^p
\end{eqnarray*}

\newpage
\begin{eqnarray*}
& & +\sum_{j=1}^{j_2} \left| a_{i2m+(i+1)+j}-\dfrac{1}{(2m)^{1/p}}\,
a_{i2m+(i+1)} \right|^p +\\
& & \left. \dfrac{2m-j_2}{2m}\left| a_{k_2\, 2m+(k_2+1)} \right|^p
\right]^{1/p} \\
& & \leq b^2 \left\| \sum_{i=1}^s\, a_i\, y_i \right\|.
\end{eqnarray*}

So the basis constant of $(y_n)$ is $b^2$ where $b$ is the basis constant
of $(x_n)$.

Since $(x_n)<(e_n^p)$, there is a constant $K\geq 1$ so that for every
sequence of scalars $(a_n)$ we have,
\begin{equation}
\left\| \sum_{n=1}^{\infty}\, a_n \, x_n \right\| \leq
K\left(\sum_{n=1}^{\infty}\, |a_n|^p\right)^{1/p}.
\label{eq:2.11}
\end{equation}

On the other hand, if $I\subset \{k 2 m+(k+1)+1, k 2 m+(k+1)+2, \ldots,
(k+1)2m+(k+1)\}$ with $|I|\leq m$,

\begin{eqnarray}
\lefteqn{\left[\sum_{j\in I}\, \left|
a_{k2m+(k+1)+j}-\dfrac{1}{(2m)^{1/p}}\, a_{k2m+(k+1)}\right|^p
+  \dfrac{2m-|I|}{2m}\,\left|a_{k2m+(k+1)}\right|^p\right]^{1/p}}\\
& & \geq \dfrac14 \left[ \sum_{j\in I}\, \left|
a_{k2m+(k+1)+j}\right|^p+\left|a_{k2m+(k+1)}\right|^p \right]^{1/p}.
\label{eq:2.12}
\end{eqnarray}
Now, if $I\subset\N$, $|I|=m$, for each $k=0, 1, \ldots$ let $I_k=I \cap
\{k2m+(k+1), k2m+(k+1)+1, \ldots, (k+1)2m+(k+1)\}$ and let $P$ be the
natural projection of $X\oplus l_p$ onto $l_p$.  From~\ref{eq:2.12} we
have,
\begin{eqnarray}
\left\| \sum_{i\in I}\, a_i \, y_i \right\| & \geq &
\left(\sum_{k=0}^{\infty}\, \left\| \sum_{i\in I_k}\, a_i\, Py_i
\right\|^p \right)^{1/p} \\
& \geq & \dfrac14 \left(\sum_{k=0}^{\infty}\, \sum_{i\in I_k}\, |a_i|^p
\right)^{1/p}.
\label{eq:2.13}
\end{eqnarray}

By~\ref{eq:2.11} and~\ref{eq:2.13}, we have that $(y_i)_{i\in I}
\approx_{8K}
(e_i^p)_{i\in I}$.

This completes the proof of the proposition.
\end{pf}

For $1<p\leq 2$, $L_p[0,1]\approx L_p[0,1]\oplus l_2$, and the Walsh
system $(W_n)$ is a basis of $L_p[0,1]$ satisfying $(W_n)<(e_n^2)$.  So
we have,

\begin{corollary}
For $1<p\leq 2$, $L_p[0,1]$ has basis property 2.
\end{corollary}

Since any space with basis property $p$ contains uniformly complemented
$l_p^n$'s, it follows that $L_1[0,1]$ fails basis property 2.  It is
immediate that every basis $(x_n)$ for a Banach space $X$ satisfies
$(e_n^0)<(x_n)<(e_n^1)$.  So for example, $L_1[0,1]$ has basis property
1.  That is, $L_1[0,1]$ has bases with uniformly bounded basis constants
every $n$--elements of which are the unit vector basis of $l_1^n$ (and
hence are well unconditional) despite the fact that $L_1[0,1]$ has no
unconditional basis.  It can be shown that $T^{(2)}$--Tsirelson's space
\cite{CS} has basis property 2 despite its not containing a subspace
isomorphic to $c_0$ or $l_p$, for any $p$.  Also, $(\sum \oplus
T^{(2)})_{l_2}$ has basis property 2.

We could strengthen this notion to ``unconditional'' basis property $p$
by requiring the unconditional basis constant to be $\leq K$ in
definition 2.9.
Although $l_p\oplus l_2$ has basis property $p$, for $2<p<\infty$ (by
proposition 2.10) it fails to have unconditional basis property $p$ by
the quantitative version of a result of Edelstein and Wojtaszczyk \cite{EW}
(see also \cite{W}, \cite{CKT}).  Also, $T^{(2)}$ fails to have
unconditional basis property 2 by the uniqueness, up to a permutation, of
the unconditional basis for $T^{(2)}$ \cite{BCLT}.  We do not know of a
non--Hilbert space with unconditional basis property 2.  Probably, such
examples exist in the class of Orlicz spaces (Probably even spaces with
``symmetric'' basis property 2).  This whole idea could warrant further
study if we could first find a good use for this concept.

Recently, V. Mascioni [Ma] has introduced an interesting variant of the
HBS problem.

\begin{problem} \cite{Ma}
If $X$ is an infinite dimensional Banach space, and every infinite
dimensional subspace of $X$ is isomorphic to its dual space, is $X$
isomorphic to a Hilbert space?
\label{pr:2.11}
\end{problem}

We could also ask the isometric version of problem~\ref{pr:2.11}.
Mascioni
[Ma] then finite dimensionalizes the problem.

\begin{defn} \cite{Ma}
A Banach space is locally self dual (LSD) if there is a constant $c$ such
that every finite dimensional subspace of $X$ is $c$--isomorphic to its
dual space.
\end{defn}

\begin{problem} \cite{Ma}
Are LSD spaces isomorphic to Hilbert spaces?
\end{problem}

We will discuss these results further in section 4.

\section{Some Positive Results}

The first major advance in this area was due to J. Bourgain \cite{B*},
who solved the finite dimensional version of the HBS problem.  Later, N.
Tomczak--Jaegermann and P. Mankiewicz \cite{MTJ1} gave the best constants
for the finite dimensional homogeneous Banach space problem.  V. D.
Milman first posed this problem in its finite dimensional form.  In the
finite dimensional setting, we cannot ask for $X$ to be isomorphic to its
subspaces.  Instead, we assume that $\dim X=n$ and for some $1\leq k\leq
n$, we assume that all $k$--dimensional subspaces of $X$ are
$K$--isomorphic.  Now we must ask for a quantitative answer: what is the
smallest constant $f(K)$ so that $X$ is $f(K)$--isomorphic to a Hilbert
space?  Clearly, $n=1$ will give no information about the Banach space
$X$.  Yet, we have a quite exact answer to this problem from \cite{B*},
and [MJT1].

\begin{thm} \cite{MTJ1}
If an $n$--dimensional Banach space $X$ has the property that all its
$[\alpha n]$--dimensional subspaces are $K$--isomorphic, for some
$0<\alpha<1$, then $X$ is $f(\alpha, K)$--isomorphic to a Hilbert space,
where $f(\alpha,K)=CK^{3/2}$, if $0<\alpha<2/3$ and $f(\alpha,K)=CK^2$,
if $2/3 <\alpha<1$, and $C=C(\alpha)$ depends on $\alpha$ only.
\end{thm}

We will come back to this theorem in a moment.  In the meantime, we
consider the only solution of a special case of the infinite dimensional
problem due to W. B. Johnson [J1].  Recall that a Banach space $X$ is
said to have the \underline{GL--property} (Gordon--Lewis property) if
every absolutely summing operator from $X$ to $L_2$ factors through
$L_1$.  Y. Gordon and D. Lewis [GL] in a landmark paper showed that every
Banach space with an unconditional basis (or even LUST) has the
GL--property.  W. B. Johnson [J1] conjectured that every Banach space has
a subspace with the GL--property.  This conjecture does not seem to have
been tested yet on the new examples of Banach spaces without
unconditional bases due to W. T. Gowers and B. Maurey
[GoM].  For a full understanding of Johnson's result, we need to recall
some definitions.

\begin{defn}
(1) A Banach space $X$ is a \underline{weak cotype $q$ space} if there is
a constant $wc_q (X)$ so that for all $n$ and all operators $u:l_2^n
\longrightarrow X$, we have $$\sup_{k\geq 1}\, k^{1/q}\, a_k\, (u) \leq
we_q (X)\, l(u),$$ where $$l(u)=\left( \int_{\Bbb{R}}\, \|u(x)\|^2\,
\gamma_n
(dx)\right)^{1/2}$$ and $\gamma_n$ is the canonical Gaussian probability
measure on $\Bbb{R}^n$, and $$a_k(u)=\inf\{\|u-s\| \left| s:l_2^n
\longrightarrow X, \mbox{ rank }\, s\leq k\right. \}.$$

(2) A Banach space $X$ is a \underline{weak type $p$ space} if there is a
constant $wt_p(X)$ so that for all $n$ and all operators
$V:X\longrightarrow l_2^n$ we have $$\sup_{k\geq 1}\, k^{1/p}\, a_k(V)
\leq wt_p(X)\, l^*(V)$$ where $$l^*(V)=\sup\{tr(uv)\left| u:l^n_2
\longrightarrow X, l(u)\leq 1\right. \}.$$

(3) A Banach space $X$ is a \underline{weak Hilbert space} if it is weak
type 2 and weak cotype 2 space.

(4) A Banach space is \underline{as--Hilbertian} (asymptotically
Hilbertian) if there is a constant $K\geq 1$ so that for every $n$, there
is a finite codimensional subspace $H$ of $X$ with the property that
every $n$--dimensional subspace of $H$ is $K$--isomorphic to $l_2^n$.
\label{de:3.2}
\end{defn}

See G. Pisier [P] for a complete treatment of these notions.  Johnson
observed;

\begin{thm} \cite{J1}
A homogeneous weak Hilbert space is isomorphic to a Hilbert space.
\label{thm:3.3}
\end{thm}

The main result from [J1] is:

\begin{thm} \cite{J1}
If $X$ and $X^*$ are homogeneous and if $X$ has the GL--property, then
$X$ is isomorphic to a Hilbert space.
\end{thm}

Returning to Theorem 3.1, we might naturally formulate an alternative
infinite dimensional HBS problem as:

If $X$ is $K$--isomorphic to all of its finite codimensional subspaces,
is $X$ $f(K)$--isomorphic to a Hilbert space?

It is easily seen that
this property is so strong that even $c_0$, $l_p$, $1\leq p<\infty$, fail
it.  In particular, a space $X$ with this property also satisfies the
property that all finite dimensional subspaces of $X$ are $K$--isomorphic
to $K$--complemented subspaces of $X$.  Also, the proof of
Theorem 3.4 shows that whenever $X$ and $X^*$ have this property
and $X$ has the GL--property, then $X$ is isomorphic to a Hilbert space.
However, this property
does not characterize a Hilbert space as our next proposition shows.

\begin{proposition}\hfill\break
\begin{enumerate}
\item If $X$ is a Banach space with a separable dual $X^*$, and
$\varepsilon
>0$, then there is a Banach space $Y$ with a separable dual containing
$X$ and $Y$ is $1+\varepsilon$--isomorphic to every subspace of $Y$ of
finite codimension.

\item If in (1), $X^{**}$ is separable, we may construct our $Y$
satisfying
(1) and so that $Y^{*}$ is also $1+\varepsilon$--isomorphic to every
subspace of $Y^*$ of finite codimension.
\end{enumerate}
\end{proposition}

\begin{pf}
Let $X_1=\left(\sum_{k=1}^{\infty}\oplus X\right)_{l_2}$ and choose
$f_n^1\in X_1^*$ so that $\{f_n^1\}$ is dense in $X^*_1$ and each $f_n^1$
appears infinitely many times in the sequence.  Let $J$ be the family of
all finite subsets of $\N$.  For each $A\in N$ let $Z(A)=\span_{n\in A}\,
f_n^1$ and let
\begin{enumerate}\itemsep 5mm
\item $X_2=X_1 \oplus_2\, \left(\sum_{A\in J}\oplus
Z(A)_{\perp}\right)_{l_2}$

It follows that

\item $X_2 \approx_{(1)} X_1 \oplus X_2,$

and if $A\in J$, then

\item $X_2 \approx_{(1)} X_2 \oplus_2 Z(A)_{\perp}$.

Now choose
$(f_n^2)$ dense in $X_2^*$ with each $f_n^2$ appearing infinitely many
times in the sequence and repeat the construction to get $X_3$.  By
induction, we construct $X_1, X_2, \ldots$.  Now let

\item $Y=\left(\sum_{n=1}^\infty \oplus X_n\right)_{l_2}$.

If follows that

\item $Y\approx_{(1)} \left(\sum_{k=n}^\infty \oplus X_k \right)_{l_2}$,
for every $n$, and

\item if $H$ is a subspace of $\left(\sum_{i=1}^n\oplus X_i
\right)_{l_2}$ of finite codimension, and $H=\left(\span_{k\in A}\,
f_k^n\right)_{\perp}$, for some $A\in J$, then
$$\left(\sum_{k=n+1}^\infty \oplus X_k \right)_{l_2} \approx_{(1)}
\left(\sum_{k=n+1}^\infty \oplus X_k \right)_{l_2} \oplus_2 H.$$

Now let $H\subset Y$ be a subspace of finite codimension and
$H^{\perp}=\span_{1\leq i\leq n}\, f_i$, where $f_i=\sum_{j=1}^\infty\,
f_{ij}$, $f_{ij}\in X_j^*$ and the basis constant of $\{f_i\}_{i=1}^n$ is
$\leq n^2$.  Choose $\varepsilon_0 < 1+\varepsilon / 2\sqrt{n}n^2$ and an
$m$ so that

\item $\|f_i-\sum_{j=1}^m\, f_{ij} \| <\varepsilon_0$, for all $i=1, 2,
\ldots, n$.

Let $g_i=\sum_{j=1}^m\, f_{ij}$, $1\leq i\leq n$ and choose ${f_k}_i^m
\in \left(\sum_{j=1}^m \oplus X_j^* \right)_{l_2}$ so that

\item $\|{f_k}_i^m -g_i\|<\varepsilon_0$.

It follows that

\item $\|{f_k}_i^m - f_i \| <2 \varepsilon_0$.

Choose a $w^*$--continuous projection $P^*$ of $Y^*$ onto $\span_{1\leq i
\leq n}\, {f_k}_i^m$ with $\|P^*\|<2 \sqrt n$.  Finally, define $T:Y^*
\longrightarrow Y^*$ by: Given $f\in Y^*$, let $P^*(f)=\sum_{i=1}^n\,
a_i\, {f_k}_i^m$ and define

\item $T(f)=\sum_{i=1}^n\, a_i\, f_i +(I-P^*)(f)$.

We now have, for all $f\in Y^*$,

\item $\|(I-T)(f)\|\leq \sum_{i=1}^n\, |a_i| \|f_i-{f_K}_i^m \| < 2
\sqrt{n} n^2 \varepsilon_0 < 1+ \varepsilon$.

So $T$ is a $w^*$
continuous isomorphism on $Y^*$ with $\|T\|\, \|T^{-1}\|\leq
(1+\varepsilon)^2$.  Let $S:T \longrightarrow Y$ be the isomorphism
satisfying $T=S^*$.  Finally, $T\left(\span_{i\leq i \leq n}\,
{f_k}_i^n\right)=H^{\perp}$ implies $S(H)\subset \left(\span_{1\leq i
\leq n} \, {f_k}_i^m \right)_{\perp}$, and dimension considerations
yields that $S$ is a $(1+\varepsilon)^2$ isomorphism of $H$ onto
$\left(\span_{1\leq i \leq n}\, {f_k}_i^n \right)_{\perp}$.  Now,
$H\approx_{(1+\varepsilon)^2} \, \left(\span_{1\leq i \leq n}\, {f_k}_i^n
\right)_{\perp}$ and so $H \approx_{(1+\varepsilon)^2}\, H_1 \oplus
\left(\sum_{j=m+1}^\infty \oplus X_j \right)_{l_2}$ where $H_1$ has
codimension $n$ in $\left(\sum_{j=1}^n \oplus X_j \right)_{l_2}$.  So by
(5) and (6), $$H\approx_{(1+\varepsilon)^2} \left(\sum_{j=m+1}^\infty
\oplus X_j\right)_{l_2} \approx_{(1)} Y.$$

For part (2) of the proposition, we alternate the above construction
between $X_n$ and $X_n^*$.
\end{enumerate}
\end{pf}

A variation of this construction carried out transfinitely would yield an
alternate proof of a result in the literature.  (This result definitely
exists in the literature despite our inability to locate it at this
time).  One
can view this as a counterexample to the ``non--separable'' HBS problem.

\begin{thm} []
There is a non--separable Banach space $X$ which is not isomorphic to a
Hilbert space, but such that $X$ is isometric to every subspace of the
same density as $X$.
\end{thm}

\section{Some Stronger Problems}

At this time, weak Hilbert space theory seems to be closing in on the HBS
problem.  Several natural questions from that area would yield a positive
solution to the HBS problem.  Since every subspace of a homogeneous
Banach space has a basis, a positive answer to the following
question,
combined with Theorem \ref{thm:3.3}, would yield a positive answer to
the HBS problem:

\begin{question}
If every subspace of $X$ has a basis, is $X$ a weak Hilbert space?
\label{q:4.1}
\end{question}

The converse of \ref{q:4.1} is also an open problem.  That is; Does
every separable weak Hilbert space have a basis?  A result of B. Maurey
and G. Pisier (which appears for the first time in [Ma]) states
that a separable weak Hilbert
space has a finite dimensional decomposition.  R. Komorowski [K] has
constructed the first weak Hilbert spaces which fail to have
unconditional bases (and they are even unconditional sums of two
dimensional subspaces).

Johnson's result \ref{thm:3.3} actually asserts that a homogeneous
Banach space which is as--Hilbertian is isomorphic to a Hilbert space.
And an earlier result of Johnson (see [P]) states that every weak Hilbert
space is as--Hilbertian.  So a weaker formulation of \ref{q:4.1} would
be:

\begin{question}
If every subspace of $X$ has the approximation property, is $X$
as--Hilbertian?
\label{q:4.2}
\end{question}

It is easily seen that there are as--Hilbertian spaces which fail the
approximation property.  Also, Johnson [J2] has exhibited a class of
Banach spaces which are as--Hilbertian, every subspace has the
approximation property (even the finite dimensional decomposition
property) but they have subspaces without bases.  There is a possible
counterexample to question~\ref{q:4.2}.  That is, the symmetric
convexified Tsirelson space [CS].  This is a Banach space with a
symmetric basis for which all $n$--dimensional subspaces are within a
fixed iterated logarithm of $l_2^n$.  It is possible, however, that every
non Hilbert space with a symmetric basis has a subspace which fails the
approximation property.

N. Tomczak--Jaegermann and P. Mankiewicz [MTJ1] obtained the following
corollary while proving some deep results in the local theory of Banach
spaces:

\begin{thm} \cite{MTJ1}
Let $X$ be a Banach space for which there exists a constant $K$ such that
every finite dimensional subspace $F$ of $X$ satisfies $bc(F)\leq K$.
Then $X$ is of weak cotype 2.
\label{thm:4.3}
\end{thm}

This Theorem gives added importance to the following question which has
been around for quite some time.

\begin{question}
Are the following equivalent for a Banach space $X$?

\begin{enumerate}
\item Every subspace of $X$ has a basis.

\item There is a constant $K\geq 1$ so that every finite dimensional
subspace of $X$ has a basis with basis constant $\leq K$.
\end{enumerate}
\label{q:4.4}
\end{question}

The importance of question~\ref{q:4.4} is that a positive solution to
$(1)\Longrightarrow (2)$ would yield that every homogeneous Banach space
is a weak cotype 2 space.  This might then be combined with a positive
solution to question~\ref{q:2.3} to solve the whole problem.

The argument of Johnson [J2] is ``local'' and shows that convexified
Tsirelson's space has both properties (1) and (2) of
question~\ref{q:4.4}.  N. Nielsen and N. Tomczak--Jaegermann [NTJ] have
shown that very weak Hilbert space with LUST also has these properties.

Recently, P. Mankiewicz and N. Tomczak--Jaegermann [MTJ3] made an
important step in resolving the above questions:

\begin{thm} \cite{MTJ3}
If a Banach space $X$ has the property that every subspace of every
quotient space of $l_2(X)$ has a Schauder basis, then $X$ is isomorphic
to a Hilbert space.
\label{thm:4.5}
\end{thm}

Neither implication in question~\ref{q:4.4} is known, which points out a
serious gap in the available techniques in Banach space theory.  Namely,
we have no reasonable way of passing results back and forth between local
theory and infinite dimensional theory.  The paper [MTJ3] should be read
not only for the main result, but also because it is the first serious
integration of local theory and infinite dimensional theory.

Mascioni [Ma] has proved the corresponding result for LSD--spaces.

\begin{thm} \cite{Ma}
If $l_2(X)$ is LSD, then $X$ is isomorphic to a Hilbert space.
\end{thm}

Another possible counterexample to question~\ref{q:4.2} is $X=\left(\sum
\oplus T^{(2)}\right)_{l_2}$, where $T^{(2)}$ is convexified Tsirelson
space (see [CS]).  This space is of type 2 but fails to be weak cotype 2,
while still satisfying: $$\inf \{q | X \mbox{ is cotype }\, q \} =2.$$
Results of N. Tomczak--Jaegermann and P. Mankiewicz [MTJ1] show that $X$
fails the finite dimensional basis property (i.e. (2) of
question~\ref{q:4.4}).  Also, results from [MTJ3] show that this space
has a subspace of a quotient space which fails to have a basis.  It is
possible, however, that every subspace of $X$ has the approximation
property.

Finally, a positive answer to the following question (plus its dual
formulation for weak cotype) would yield a positive solution to the HBS
problem:

\begin{question}
If $\sup\{p|X \mbox{ is Type } p\}=p_0$, does $X$ contain a subspace of
weak type $p_0$?
\label{q:4.5}
\end{question}

\section{Some Basic Questions}

Let us return to some of the open questions of section 2.  It seems
strange that we do not know if a homogeneous Banach space is reflexive,
especially in light of formula (1) of section 2.  R. C. James [Ja] has
given us a class of Banach spaces with type which fail to be reflexive.
But a positive answer to the following question would show
that homogeneous Banach spaces are reflexive.

\begin{question}
If a Banach space $X$ is of type $p$, for some $1<p\leq 2$, must $X$
contain a reflexive subspace?
\label{q:5.1}
\end{question}

If a Banach space has a subspace with LUST, the answer is ``yes'' (see
[LT2]).  Actually, a stronger conclusion could hold:

\begin{question}
Does every Banach space of type $p$ for some $1<p\leq 2$, contain a
superreflexive subspace?
\label{q:5.2}
\end{question}

The cotype version of question~\ref{q:5.1} also seems to be unasked and
unanswered:

\begin{question}
If $X$ is of cotype $q$, for some $2 \leq q \leq \infty$, must $X$
contain a separable dual space?
\label{q:5.3}
\end{question}

An equivalent formulation of question~\ref{q:5.3} would be to ask if $X$
must contain a boundedly complete basic sequence.  These questions are special cases of the question of H. P. Rosenthal:
Does every Banach space contain either a reflexive subspace or a subspace
isomorphic to $c_0$ or $l_1$?  In fact, question~\ref{q:5.2} is the
``uniform'' or ``local'' version of Rosenthal's question.  That is,
question~\ref{q:5.2} is equivalent to:

{\bf Question 5.5.2'}
{\it Does every Banach space contain either a superreflexive subspace or
subspaces uniformly isomorphic to $l_1^n$ or $l_{\infty}^n$, for every
$n=1,2, \ldots$?}

In this setting, the ``$l_{\infty}^n$'' in the question is redundant.

Recently W. T. Gowers [Go] has constructed a Banach space not containing
$c_0$, $l_1$ or any reflexive subspace.  In fact, Gowers' space has no
subspace with a separable dual space.  It is possible that refinements of
this example will give counterexamples to the above problems.

\section{A Construction}

What is the ``best'' subspace we can find inside of every Banach space?
The major problem here is not just to answer the question, but to
formulate the question.  We now know that almost every Banach space
contains a subspace which fails the approximation property (see [LT2],
Chapter 1g).  Thanks to W. B. Johnson [J2] we also know of non--Hilbert
spaces for which every subspace of every quotient space has a basis.  We
also know that every Banach space contains a basic sequence but may not
contain an unconditional basic sequence [GoM].  The best subspace we
could hope to find in a Banach space is a subspace isomorphic to $c_0$,
or $l_p$, $1 \leq p <\infty$.  Not just because we understand these
spaces better than any others, but because they have the property that
they embed (complementably) into every infinite dimensional subspace of
themselves.  That is, we can recover all of the properties of the whole
space inside of every subspace.  But in 1972, B. S. Tsirelson [T] (see
also [CS]) showed that there are Banach spaces which do not contain
copies of $c_0$ or $l_p$, $1\leq p <\infty$.  This example quickly
blossomed into an ``industry'' [CS] and even today has its place in the
recent exciting solutions to the unconditional basic sequence problem
[GoM], the $c_0$ $l_1$, reflexive space counterexample [Go], and the
distortion problem
[S1],
[OS],
[MiTJ]. H. P.
Rosenthal then raised the question whether every Banach space $X$ might
contain a subspace $Y$ which embeds into every one of its infinite
dimensional subspaces?  Such a space $Y$ is called \underline{minimal}.
This was the ``best'' subspace we could hope for at the time.  In 1982,
P. G. Casazza and E. Odell [C0] (see also [CS]) showed that Tsirelson's
space contains no minimal subspaces.  As of this writing, we know of only
two new classes of minimal Banach spaces (besides subspaces of $c_0$,
$l_p$, $1\leq p < \infty$) [CJT], [S2], and only the second is
complementably minimal.  So, we now know that we cannot find, in every
Banach space $X$, a subspace $Y$ so that infinite dimensional properties
of $Y$ are invariant under passing to further subspaces.

Our next approach would be to look for a ``locally best'' subspace in
every Banach space, i.e. a subspace $Y$ of $X$ so that $Y$ is crudely
finitely representable in every one of its subspaces.  This is not
possible either as Tsirelson's space again fails this property by the
argument of [OS].  Since Tsirelson's space seems to be blocking all our
efforts, let's see what property this space does have.  The Tsirelson
space $T_p$ enjoys the property that there is a constant $K$ so that for
every $n$ there is a subspace $H$ of $T_p$ of condimension $n$ and every
$n$--dimensional subspace of $H$, $K$--embeds into every infinite
dimensional subspace of $T_p$.  From our earlier examples, this is the
``best'' we can hope for in an arbitrary Banach space.  Our next theorem
states that, indeed, every Banach space does contain such a subspace (and
with $K$ arbitrarily close to 1).

Recall that a Banach space $Y$ (finite or infinite dimensional)
\underline{almost isometrically}\hfill\break \underline{embeds} into a
Banach space
$X$ if for
every $\varepsilon >0$, there is a subspace $Z\subset X$ and an operator
$T:Y \longrightarrow Z$ so that $\|T\|\, \|T^{-1}\| \leq 1+\varepsilon$.
We can now state the theorem.

\begin{thm}
For every Banach space $X$, for every $\varepsilon_n \downarrow 0$ and
for every $f:N\longrightarrow N$ there is a subspace $Y\subset X$ which
the following properties:

\begin{enumerate}
\item $Y$ has a normalized basis $(y_n)$ with basis constant $\leq 1 +
\varepsilon_0$,

\item Every $E\subset \span_{n\leq i < \infty}\,  y_i$, with $\dim E \leq
f(n)$, is $1+\varepsilon_n$--isomorphic to a Banach space $F$, and $F$
almost isometrically embeds into every infinite dimensional subspace of
$Y$,

\item Every block basis $(z_i)_{i=1}^{f(n)}$ of $(y_i)_{i=n}^{\infty}$ is
$1+\varepsilon_n$--equivalent to a basis $\{w_i\}_{i=1}^{f(n)}$ of a
Banach space $F$, and $\{w_i\}_{i=1}^{f(n)}$ is almost isometrically
equivalent to a block basis of every basic sequence in $Y$.
\end{enumerate}
\label{thm:6.1}
\end{thm}

Before we prove Theorem~\ref{thm:6.1}, let us consider some of its
consequences.  We have immediately from Definition~\ref{de:3.2} and
Theorem~\ref{thm:6.1}:

\begin{corollary}
If $X$ is an infinite dimensional Banach space, there is an infinite
dimensional subspace $Y$ of $X$ so that

\begin{enumerate}
\item $wt_p(Z)=wt_p(Y)$,

\item $wc_q(Z)=wc_q(Y)$,
\end{enumerate}

for every infinite dimensional subspace $Z$ of $Y$ and every $1\leq p
\leq 2$, $2 \leq q < \infty$.
\label{co:6.2}
\end{corollary}

It follows that $Y$ has type (respectively, cotype) if and only if $Y$
has a subspace with type (respectively, cotype).  Also, $Y$ is a weak
Hilbert space if and only if $Y$ contains a weak Hilbert space.

\begin{corollary}
For every Banach space $X$ there is a subspace $Y\subset X$ so that every
spreading model of $Y$ (built from a weakly null sequence)
is finitely representable in every
(other) infinite dimensional subspace of $Y$.
\label{co:6.3}
\end{corollary}

\begin{pf}
This is immediate since spreading models of subspaces of $Y$
are finitely representable in $\span_{n\leq i < \infty}\,
y_i$, for every $n=1,2,\ldots$.
\end{pf}

Thus, if $Y$ has $c_0$ as a spreading model (respectively $l_1$) then $Y$
has no subspaces with cotype (respectively, type).
Recall that a basic sequence $(x_n)$ is \underline{block finitely}
\underline{representable} in a basic sequence $(y_n)$ if for every $n$
and every
$\varepsilon > 0$ there is a block basis $(z_i)_{i=1}^n$ of $(y_i)$ which
is $1+\varepsilon$--equivalent to $(x_i)_{i=1}^n$.  In this notation,
Krivine's Theorem [Kr] (see also [MS]) says that for every basic sequence
$(x_n)$ in a Banach space $X$ there is a $1\leq p < \infty$ (or for
$c_0$) so that the unit vector basis of $l_p$ (or $c_0$) is block
finitely representable in $(x_n)$.  Property (3) of Theorem~\ref{thm:6.1}
implies that whenever the unit vector basis of $l_p$ is block finitely
representable on $(y_n)$ then it is block finitely representable on every
basic sequence of $Y$.  This result was first observed by H. P. Rosenthal
[R].

\begin{corollary}
For every Banach space $X$, there is a subspace $Y$ and a $1 \leq p <
\infty$ (or $c_0$) so that the unit vector basis of $l_p$ (or $c_0$) is
block finitely representable in every basic sequence in $Y$.  Moreover,
if $l_p$ (or $c_0$) is block finitely representable in one basic sequence
in $Y$, then it is block finitely representable in every basic sequence
in $Y$.
\label{co:6.4}
\end{corollary}

Recall that a Banach space $X$ is \underline{$K$--crudely finitely
representable} in a Banach space $Y$ if for every finite dimensional
subspace $E\subset X$ there is a subspace $F\subset Y$ with $d(E, F)\leq
K$.  If for every $\varepsilon >0$ and every finite dimensional subspace
$E\subset X$ there is a subspace $F\subset Y$ with $d(E, F)\leq 1
+\varepsilon$, we say $X$ is \underline{finitely representable} in $Y$.

\begin{corollary}
Suppose the separable Banach space $X$ is $K$--crudely finitely
representable in every infinite dimensional subspace of $X$.  Then there
is an equivalent norm $1\cdot1$ on $X$ so that $(X, 1 \cdot 1)$ is
finitely representable in every infinite dimensional subspace of $X$.
\label{co:6.5}
\end{corollary}

\begin{pf}
Choose the subspace $Y$ of $X$ with basis $(y_i)$ from
Theorem~\ref{thm:6.1}.  Choose finite dimensional subspaces $E_1 \subset
E_2 \subset \cdots$ whose union is dense in $X$.  For each $i, j=1,2,
\ldots$, choose $F_{ij}\subset \span_{l_(j)\leq k < \infty}\, y_{k}$,
where $l (j) \geq \dim F_{ij}$,  and $d(E_i, F_{ij})\leq K$.  By
switching to a subsequence, we may assume $\lim_{j \to \infty}\,
F_{ij}=F_{i}$, for every $j=1,2, \ldots$ (the limit in Banach--Mazur
distance) and $F_i$ is finitely representable in every infinite
dimensional subspace of $X$, by Theorem~\ref{thm:6.1} (2).  For each
$i=1,2, \ldots$ let $T_i:E_i \longrightarrow F_i$ be a $K$--isomorphism.
By switching to a subsequence again, we may assume that for every $x\in
X$, $|x|=\lim_{i \to \infty} \|T_i\, x\|$ exists.  Then $1 \cdot 1$ is an
equivalent norm on $X$ and clearly $(X, 1 \cdot 1)$ is finitely
representable in $X$ (again by Theorem~\ref{thm:6.1}).
\end{pf}

If the above corollary had an infinite dimensional analogue, it could be
quite useful for working in minimal Banach spaces.  But our proof is
local and we have not found a generalization of it for this case.  Also,
it would be much better if we could show that $(X, 1 \cdot 1)$ is
finitely representable in every infinite dimensional subspace of itself.

Again, we might hope that such a strong property would characterize a
Hilbert space.  But, it is easily seen that $c_0$ has the property that
every finite dimensional subspace \underline{isometrically} embeds into
every infinite dimensional subspace.

We are now ready for the proof of Theorem~\ref{thm:6.1}.

{\sll Proof of Theorem~\ref{thm:6.1}}.\,\,
Fix $\varepsilon > 0$ and $n\in N$.  Choose $n$--dimensional Banach
spaces $H_1, H_2, \ldots, H_m$ which are $1+\varepsilon$--dense in the
set of all $n$--dimensional Banach spaces in the Banach--Mazur distance.
We now ask: $(+)$ Is $H_1$, $1+\varepsilon$--embeddable into every
infinite dimensional subspace of $X$?

If the answer is ``yes'', put $H_1$ in the set $A$ and go to $H_2$.  If
the answer is ``no'', put $H_1$ in the set $B$ and replace $X$ by an
infinite dimensional subspace $X_1$ of $X$ so that $X_1$ has no subspace
$1+\varepsilon$--isomorphic to $H_1$.  Now go to $H_2$ and start over.
After $m$--steps, we are left with a partition of the $H_i$'s into two
groups, $A=\{H_i\}_{i\in I}$ and $B=\{H_i\}_{i\in J}$, $I\wedge J=\phi$,
and an infinite dimensional subspace $Z$ of $X$ so that every $H_i$, for
$i\in I$, is $1+\varepsilon$--embeddable into every infinite dimensional
subspace of $Z$ while no $H_i$, $i\in J$, is $1+\varepsilon$--embeddable
into $Z$.

It follows that for any $E\subset Z$, $\dim E=n$, $E$ is
$1+\varepsilon$--isomorphic to $H_i$, for some $i\in I$, and hence $E$ is
$(1+\varepsilon)^2$--embeddable into every infinite dimensional subspace
of $Z$.  Observe that this property is maintained if we switch to any
infinite dimensional subspace of $Z$.  Hence, given $f:N \longrightarrow
N$ and $\varepsilon_n \downarrow 0$, we can inductively carry out this
procedure to obtain infinite dimensional subspaces of $X$, $X\supset Z_1
\supset Z_2 \supset \cdots$ so that every $f(n)$--dimensional subspace of
$Z_n$ is $1+\varepsilon_n$--embeddable into every infinite dimensional
subspace of $Z_n$.  Now choose $y_n \in Z_n$ so that $(y_n)$ is a
$1+\varepsilon_0$--basic sequence in $X$.  Now, if $E \subset \span_{n
\leq
i}\, y_i$, then there are subspaces $F_k \subset \span_{n+k \leq i}\,
y_i$ so that $d(E, F_k) \leq 1+\varepsilon_n$.  By switching to a
subsequence, we may assume $\lim_{k \to \infty}\, F_k=F$ (in
Banach--Mazur distance).  So $d(E, F)\leq 1 +\varepsilon_n$ but $F_k$ is
$1+\varepsilon_{n+k}$ embeddable into every infinite dimensional subspace
of $Y$.  Hence, $F$ is finitely representable in every infinite
dimensional subspace of $Y$.  This concludes the construction for part
(1) of Theorem~\ref{thm:6.1}.  Again, note that part (1) of the Theorem
holds if $(y_n)$ is replaced by any block basis of $(y_n)$.

Part (2) of Theorem~\ref{thm:6.1} is proved in a similar manner.
Again, fixing $n$ and $\varepsilon >0$, we list out
$(x_{1i})_{i=1}^n, (x_{2i})_{i=1}^n, \ldots, (x_{mi})_{i=1}^n$ with the
property: Every normalized basis $(z_i)_{i=1}^n$, with basis constant
$\leq 2$, for a Banach space $X$, is $1+\varepsilon$--equivalent to one
of the $(x_{ki})_{i=1}^n$.  Using our basis $(y_i)$ from part (1), we
ask: Is $(x_{1i})_{i=1}^n$, $1+\varepsilon$--equivalent to a block basis
of every infinite block basis of $Y$?

As before, if the answer is ``yes'', put $(x_{1i})$ in the set $A$ and go
on to $(x_{2i})$.  If the answer is ``no'', put $(x_{1i})$ in the set $B$
and replace $(y_i)$ by a block basis $(y_i^1)$ of $(y_i)$ so that
$(x_{1i})$ is not $1+\varepsilon$--equivalent to any block basis of
$(y_i^1)$.  Now go to $(x_{2i})$ and start over.  After $m$ steps, we
arrive at a block basis $(z_i)$ of $(y_i)$ so that every
$(x_{ji})_{i=1}^n$ in $A$ is $1+\varepsilon$--equivalent to a block basis
of every block basis of $(z_i)$ and $(x_{ji})_{i=1}^n$ in $B$ is not
$1+\varepsilon$--equivalent to any block basis of $(z_i)$.  As in part
(1), we perform this construction inductively to produce successive
block bases of block bases $(z_i^1), (z_i^2), \ldots$ with the above
property for $1+\varepsilon_n$.  Then $(z_n^n)$ is the required block
basis.  Relabeling $(z_n^n)$ as $(y_n)$, we see that we now have property
(3) of Theorem~\ref{thm:6.1} while maintaining property (2).  This
completes the proof of the Theorem.\qquad $\square$

Theorem~\ref{thm:6.1} can be partially strengthened in several
directions.  For one, if we let $Y_n=\span_{n\leq i}\, y_i$, then for
every $E\subset Y_n^*$, $\dim E\leq f(n)$, and for every $m$, there is an
$F\subset Y_m^*$ with $d(E, F) \leq 1+\varepsilon$.  It is easily seen
using $X=l_1$ that we cannot expect to $1+\varepsilon$ embed $E$ into
every infinite dimensional subspace of $Y^*$.  With significantly more
effort, we can show that our blocks in Theorem~\ref{thm:6.1} are
constructable in a very strong way.  That is, whenever $(z_i)_{i=1}^n$ is
a block basis of $(y_i)_{i=n}^\infty$, then for every $m_1$, there is a
$w_1 \in \span_{m_1 \leq k} \, y_k$, so that for every $m_2$, there is a
$w_2 \in \span_{m_2 \leq k}\, y_k, \ldots$, so that for every $m_n$
there is a $w_n \in \span_{m_n \leq k}\, y_k$ and $(w_i)_{i=1}^n$ is
$1+\varepsilon_n$--equivalent to $(z_i)_{i=1}^n$.

Now, a final word about Krivine's Theorem.  Until recently, there was the
possibility for strengthening this powerful and useful result.  First,
given $\varepsilon > 0$, and $n\in N$, we might look in every Banach
space $X$ for a $p$ and a basic sequence $(x_k)$ in $X$ so that every
block basis $(y_i)_{i=1}^n$ of $(x_k)$ is $1+\varepsilon$--equivalent to
the unit vector basis of $l_p^n$.  It can be shown that Tsirelson's space
$T^{(q)}$, fails this property for every $1\leq q < \infty$.  But
Tsirelson's space has this property for $\varepsilon = 1/2$.
Unfortunately, Schlumprecht's space [S2] fails this property for every
$\varepsilon > 0$.  Our final hope for a strengthening of Krivine's
Theorem has just fallen to the Gowers'--Maurey spaces [GoM].  That is, in
general we
cannot get basic sequences in a Banach space $X$ for which small numbers
of blocks are even well unconditional.  In particular, carrying out the
construction of proposition~\ref{thm:6.1} in the ``unconditional''
Gowers' spaces [Go2] we obtain,

\begin{proposition}
For every $K\geq 1$, there is an $n$ and $n$--vectors $(x_i)_{i=1}^n$ in
a Banach space $X$ and an infinite dimensional Banach space $Y$ with an
unconditional basis so that $(x_i)$ is 1--block finitely representable on
every basic sequence in $Y$, yet $(x_i)$ is not $K$--unconditional.
\label{p:6.6}
\end{proposition}

%%%%%%%%%%%%%%%%%%%%%%%%%%%%%%%%%%%%%%%%%%%%%%%%%%%%%%%%%%%%%%%%%%%%%%%%
%%  When BibTeX is used to create a bibliography, the LaTeX input file
%%  must include three commands of the following form.  When run through
%%  LaTeX, they will create the .aux file required by BibTeX to select
%%  the appropriate items from the specified .bib file and create a
%%  .bbl file.
%%
%\nocite{*}
%\bibliographystyle{amsplain}
%\bibliography{amsl-bib}
%%
%%  After the article is completed and the bibliography checked, to make
%%  sure that no changes need to be made in the original .bib file which
%%  require BibTeX to be rerun, the \nocite, \bibliographystyle and
%%  \bibliography commands should be commented out or removed, and the
%%  contents of the .bbl file inserted in their place, as done here.

%%%%%%%%%%%%%%%%%%%%%%%%%%%%%%%%%%%%%%%%%%%%%%%%%%%%%%%%%%%%%%%%%%%%%%%%
%%  The contents of the file amsl-art.bbl follow:

%\makeatletter \renewcommand{\@biblabel}[1]{\hfill#1.}\makeatother
%\newcommand{\bysame}{\leavevmode\hbox to3em{\hrulefill}\,}

%%  End of included file amsl-art.bbl.
%%%%%%%%%%%%%%%%%%%%%%%%%%%%%%%%%%%%%%%%%%%%%%%%%%%%%%%%%%%%%%%%%%%%%%%%

\end{document}